\documentclass[12pt,notitlepage]{article}
\usepackage{amssymb}
\usepackage{amsmath}
\usepackage{amsfonts}
\usepackage{mitpress}
\usepackage{hyperref}

\newdimen\dummy
\dummy=\oddsidemargin
\input{tcilatex}
\begin{document}

\title{Numerical investigation of the solutions of Schr\"{o}dinger equation
with exponential cubic B-spline finite element method}
\author{Ozlem Ersoy$^{a}$, Idris Dag$^{a}$ and Ali Sahin$^{b}$ \\
$^{a}${\small Department of Mathematics and Computer, Eskisehir Osmangazi
University, 26480, Eskisehir, Turkey.}\\
$^{b}${\small Department of Mathematics, Aksaray University, 68100, Aksaray,
Turkey.}}
\maketitle

\begin{abstract}
In this paper, we investigate the numerical solutions of the cubic nonlinear
Schr\"{o}dinger equation via the exponential B-spline collocation method.
Crank-Nicolson formulas are used for time discretization of the target
equation. A linearization technique is also employed for the numerical
purpose. Four numerical examples related to single soliton, collision of two
solitons that move in opposite directions, the birht of standing and mobile
solitons and bound state solution are considered as the test problems. The
accuracy and the efficiency of the purposed method are measured by L$%
_{\infty }$ error norm and conserved constants. The obtained results are
compared with the possible analytical values and those in some earlier
studies.
\end{abstract}



\section{Introduction}

One of the most interesting universal equation in physical studies is the
Schr\"{o}dinger equation that describes the quantum state of a physical
system. Since the equation used in quantum mechanics is too general, there
are some different versions of the Schr\"{o}dinger equation in scientific
studies for the modelling several phsical phenomena such as the propagation
of optical pulses, superconductivity, waves in water and plasmas and self
focusing in laser pulses. Here we focus on one of the specific form of the
Schr\"{o}dinger equation known as the time dependent cubic nonlinear Schr%
\"{o}dinger equation (NLS) which describes the optical pulse propagation in
optical fibers.

The cubic NLS equation is given in one dimension as follows: 
\begin{equation}
iU_{t}+U_{xx}+q\left\vert U\right\vert ^{2}U=0,\text{ \ \ \ \ }-\infty
<x<\infty ,\text{ \ }t>0  \label{NLS}
\end{equation}%
where $i$ is the imaginary unit, $q$ is the parameter for the self phase
modulation, $U$ is a complex-valued function which shows the evolution of
slowly varying wave train in a stable dispersive physical system with no
dissipation and $U_{t}$ is the amplitude of the pulse envelope. To complete
the usual classical mathematical statement of the problem, the initial and
the boundary conditions are chosen as to be%
\begin{equation}
U\left( x,0\right) =f\left( x\right) ,\text{ \ \ \ \ }-\infty <x<\infty
\label{IC}
\end{equation}%
and%
\begin{equation}
\lim_{x\rightarrow \pm \infty }U\left( x,t\right) =0,\text{ \ \ \ \ }t\geq 0.
\label{BC}
\end{equation}

There are many analytical and numerical studies on Eq.(\ref{NLS}) in the
literature. Different kinds of numerical techniques such as finite
difference (\cite{Twizell}, \cite{Brastsos}), finite element (\cite{I1999}, 
\cite{LRT1993}, \cite{Saka}, \cite{idag}, \cite{gar}) and Adomian
decomposition (\cite{Bratsos}) methods have been applied to Schr\"{o}dinger
equation in these studies. In last decade, variational iteration\cite%
{Sweilam}, differential quadrature\cite{Korkmaz}, cubic non-polynomial spline%
\cite{Danaf}, parametric cubic spline\cite{BinLin} and time-splitting
pseudo-spectral domain decomposition\cite{Mehdi} methods have been presented
for numerical solutions of the cubic NLS equation.

Spline functions and numerical methods where splines are used for numerical
approximation are also well studied area in applied mathematics. The first
mathematical reference to splines is the early work of Schoenberg\cite%
{Schoenberg} who revealed that splines have powerful approximation
properties. Subsequently, many approximation methods have been employed \cite%
{Ahlberg}. A spline function is a sufficiently smooth piecewise function$.$
It possesses a high degree of smoothness at the knots. A B-spline is a
special spline function that play an important role in approximation and
geometric modeling. They are used in data fitting, computer-aided design,
automated manufacturing and computer graphics. In particular, after de
Boor's \cite{deBoor} results about B-splines, spline techniques became
popular for a broad range of applications \cite{Hollig}. Most properties and
an efficient construction of B-splines can be found in \cite{deBoor}. Due to
their some attractive properties such as having compact support and yielding
numerical schemes with a high resolving power, B-splines are also widely
used in differential problems. Because of having compact support, using
B-splines in numerical solution of differential equations leads to sparse
matrix systems. The approximation of differential problems with B-splines is
obtained by the method of weighted residual, of which the Galerkin and
collocation methods are particular cases. The Galerkin method is the most
widely used method for B-spline approximations on the other hand, the
collocation method represents an economical alternative since it only
requires the evaluation at grid points \cite{Botella}. Exponential B-splines
lead to accurate\ numerical results and there are relatively less studies in
which exponential B-splines considered for the approximation. The main
objective of this paper is to construct an efficient method with the usage
of exponential cubic B-splines for the numerical investigation of cubic NLS
equation.

This paper is organized as follows: Section 2 is devoted to the numerical
method. Introducing the exponential B-splines and the application of
collocation method are given in that section. The numerical testing and the
comparisons on the examples are studied in Section 3. Finally, a conclusion
is presented in the last section.

\section{Numerical method}

Let us start with the construction of our mesh, to build on the numerical
method on it. For the computational purpose, we should restrict the solution
domain from being infinite domain to be a finite interval $[a,b].$ Since the
boundary condition (\ref{BC}) indicates that the solutions are negligibly
small outside of a finite interval, instead of physical conditions (\ref{BC}%
), we can consider the artificial boundary conditions

\begin{equation*}
\begin{tabular}{l}
$U(a,t)=U_{x}(a,t)=U_{xx}(a,t)=0,$ \\ 
\\ 
$U(b,t)=U_{x}(b,t)=U_{xx}(b,t)=0.$%
\end{tabular}%
\end{equation*}%
Then the\ uniform mesh is constructed by%
\begin{equation*}
a=x_{0}<x_{1}<\ldots <x_{N}=b
\end{equation*}%
where $x_{i}$ are knots and $h=x_{i}-x_{i-1},$ $i=1,\ldots ,N$ is the mesh
size.

\subsection{Exponential cubic B-splines}

Over the above mesh, the exponential cubic B-spline, $B_{i}${}$(x),$ is
defined by%
\begin{equation}
B_{i}(x)=\left\{ 
\begin{array}{ll}
b_{2}\left( \left( x_{i-2}-x\right) -\dfrac{1}{p}\left( \sinh (p\left(
x_{i-2}-x\right) )\right) \right) & \left[ x_{i-2},x_{i-1}\right] , \\ 
a_{1}+b_{1}\left( x_{i}-x\right) +c_{1}\exp \left( p\left( x_{i}-x\right)
\right) +d_{1}\exp \left( -p\left( x_{i}-x\right) \right) & \left[
x_{i-1},x_{i}\right] , \\ 
&  \\ 
a_{1}+b_{1}(x-x_{i})+c_{1}\exp \left( p\left( x-x_{i}\right) \right)
+d_{1}\exp \left( -p\left( x-x_{i}\right) \right) & \left[ x_{i},x_{i+1}%
\right] , \\ 
b_{2}\left( (x-x_{i+2})-\dfrac{1}{p}(\sinh \left( p\left( x-x_{i+2}\right)
\right) )\right) & \left[ x_{i+1},x_{i+2}\right] , \\ 
0 & \text{otherwise.}%
\end{array}%
\right.  \label{B_spline}
\end{equation}%
where%
\begin{equation*}
\begin{array}{l}
a_{1}=\dfrac{phc}{phc-s},\text{ }b_{1}=\dfrac{p}{2}\left( \dfrac{c(c-1)+s^{2}%
}{(phc-s)(1-c)}\right) ,\text{ }b_{2}=\dfrac{p}{2(phc-s)}, \\ 
\\ 
c_{1}=\dfrac{1}{4}\left( \dfrac{\exp (-ph)(1-c)+s(\exp (-ph)-1)}{(phc-s)(1-c)%
}\right) , \\ 
\\ 
d_{1}=\dfrac{1}{4}\left( \dfrac{\exp (ph)(c-1)+s(\exp (ph)-1)}{(phc-s)(1-c)}%
\right) .%
\end{array}%
\end{equation*}%
and $s=\sinh (ph),$ $c=\cosh (ph)$ and $p$ is a free parameter that should
be determined in computations.

A standard exponential cubic B-spline$\ $is shown in Fig.1 for $p=1.$ Each $%
B_{i}(x)$ has same shape and same size, so it is easy to see that each
exponential cubic B-spline covers four successive interval such that each
element is covered by four sequential exponential cubic B-splines.
Exponential cubic B-spline $B_{i}(x)$ and its first two derivatives are
continuous on interval $\left[ x_{i-2},x_{i+2}\right] .$

\bigskip 
\begin{equation*}
\begin{tabular}{c}
\FRAME{itbpF}{2.9257in}{2.9257in}{0in}{}{}{fig1.bmp}{\special{language
"Scientific Word";type "GRAPHIC";maintain-aspect-ratio TRUE;display
"USEDEF";valid_file "F";width 2.9257in;height 2.9257in;depth
0in;original-width 3.2197in;original-height 3.2197in;cropleft "0";croptop
"1";cropright "1";cropbottom "0";filename 'fig1.bmp';file-properties
"XNPEU";}} \\ 
\textbf{Fig.1:} Exponential cubic B-spline for $p=1.$%
\end{tabular}%
\end{equation*}

\bigskip

The nodal values and the principle two derivatives at the knots are
tabulated in Table 1. These values is going to use in application of the
numerical method.

\bigskip 
\begin{equation*}
\begin{tabular}{llclclclclc}
\multicolumn{11}{l}{\textbf{Table 1}} \\ 
\multicolumn{11}{l}{The nodal values and the principle two derivatives at
the knots} \\ \hline
&  & $x_{i-2}$ &  & $x_{i-1}$ &  & $x_{i}$ &  & $x_{i+1}$ &  & $x_{i+2}$ \\ 
\cline{3-11}
$B_{i}$ &  & $0$ &  & $\dfrac{s-ph}{2(phc-s)}$ &  & $1$ &  & $\dfrac{s-ph}{%
2(phc-s)}$ &  & $0$ \\ 
&  &  &  &  &  &  &  &  &  &  \\ 
$B_{i}^{^{\prime }}$ &  & $0$ &  & $\dfrac{p(1-c)}{2(phc-s)}$ &  & $0$ &  & $%
\dfrac{p(c-1)}{2(phc-s)}$ &  & $0$ \\ 
&  &  &  &  &  &  &  &  &  &  \\ 
$B_{i}^{^{\prime \prime }}$ &  & $0$ &  & $\dfrac{p^{2}s}{2(phc-s)}$ &  & $-%
\dfrac{p^{2}s}{phc-s}$ &  & $\dfrac{p^{2}s}{2(phc-s)}$ &  & $0$ \\ \hline
\end{tabular}%
\end{equation*}

\bigskip

\subsection{The finite element collocation method}

Considering dirac delta functions as the weighted functions in the weighted
integral of the residual leads to the situation that the residual at each
point $x_{i}$ in the domain is forced to be exactly zero. Then a system of $%
N $ residual equations is obtained in finite element collocation method.
Therefore the direct substitution of the approximation into the differential
equation is the main idea behind our numerical method.

Before the implentation of the numerical scheme, we first decompose the
governing equation (\ref{NLS}) into its real and imaginary parts as follows:%
\begin{equation}
U(x,t)=r(x,t)+is(x,t),  \label{Decomposed}
\end{equation}%
where $r(x,t)$ and $s(x,t)$ are both real valued functions. The above
decomposition yields a pair of real-valued equations, i.e.%
\begin{equation}
s_{t}-r_{xx}-q(r^{2}+s^{2})r=0,  \label{NLS_1}
\end{equation}%
\begin{equation}
r_{t}+s_{xx}+q(r^{2}+s^{2})s=0.  \label{NLS_2}
\end{equation}

Since the set of $\{B_{-1}(x),$ $B_{0}(x),\cdots ,B_{N+1}(x)\}$ forms a
basis \cite{Prenter} for the functions defined over the solution domain, the
approximations for both $r(x,t)$ and $s(x,t)$ can be constructed as%
\begin{equation*}
\begin{tabular}{lll}
$r_{N}(x,t)=\sum\limits_{i=-1}^{N+1}\delta _{i}(t)B_{i}(x),$ &  & $%
s_{N}(x,t)=\sum\limits_{i=-1}^{N+1}\phi _{i}(t)B_{i}(x)$%
\end{tabular}%
\end{equation*}%
where $\delta _{i}$ and $\phi _{i}$ are time dependent unknown parameters
that should be determined from the system of residual equations. The
derivatives of these global approximations calculated by%
\begin{equation*}
\begin{tabular}{lll}
$r_{N}^{\prime }(x,t)=\sum\limits_{i=-1}^{N+1}\delta _{i}(t)B_{i}^{\prime
}(x),$ &  & $s_{N}^{\prime }(x,t)=\sum\limits_{i=-1}^{N+1}\phi
_{i}(t)B_{i}^{\prime }(x),$ \\ 
&  &  \\ 
$r_{N}^{\prime \prime }(x,t)=\sum\limits_{i=-1}^{N+1}\delta
_{i}(t)B_{i}^{\prime \prime }(x)$ &  & $s_{N}^{\prime \prime
}(x,t)=\sum\limits_{i=-1}^{N+1}\phi _{i}(t)B_{i}^{\prime \prime }(x).$%
\end{tabular}%
\end{equation*}%
Usage of these approximations togetger with the related values in Table 1
gives the following expressions:%
\begin{equation*}
\begin{tabular}{l}
$r_{i}=r(x_{i},t)=\dfrac{s-ph}{2(phc-s)}\delta _{i-1}+\delta _{i}+\dfrac{s-ph%
}{2(phc-s)}\delta _{i+1},$ \\ 
\\ 
$r_{i}^{\prime }=r^{\prime }(x_{i},t)=\dfrac{p(1-c)}{2(phc-s)}\delta _{i-1}+%
\dfrac{p(c-1)}{2(phc-s)}\delta _{i+1}$ \\ 
\\ 
$r_{i}^{\prime \prime }=r^{\prime \prime }(x_{i},t)=\dfrac{p^{2}s}{2(phc-s)}%
\delta _{i-1}-\dfrac{p^{2}s}{phc-s}\delta _{i}+\dfrac{p^{2}s}{2(phc-s)}%
\delta _{i+1}.$%
\end{tabular}%
\end{equation*}%
and%
\begin{equation*}
\begin{tabular}{l}
$s_{i}=s(x_{i},t)=\dfrac{s-ph}{2(phc-s)}\phi _{i-1}+\phi _{i}+\dfrac{s-ph}{%
2(phc-s)}\phi _{i+1},$ \\ 
\\ 
$s_{i}^{\prime }=s^{\prime }(x_{i},t)=\dfrac{p(1-c)}{2(phc-s)}\phi _{i-1}+%
\dfrac{p(c-1)}{2(phc-s)}\phi _{i+1},$ \\ 
\\ 
$s_{i}^{\prime \prime }=s^{\prime \prime }(x_{i},t)=\dfrac{p^{2}s}{2(phc-s)}%
\phi _{i-1}-\dfrac{p^{2}s}{phc-s}\phi _{i}+\dfrac{p^{2}s}{2(phc-s)}\phi
_{i+1}.$%
\end{tabular}%
\end{equation*}%
Time discretization of decomposed system (\ref{NLS_1}) and (\ref{NLS_2}) can
be achived by Crank-Nicolson appoximation such that%
\begin{equation}
\begin{array}{r}
\dfrac{s^{n+1}-s^{n}}{\Delta t}-\dfrac{r_{xx}^{n+1}+r_{xx}^{n}}{2}+q\dfrac{%
((r^{2}+s^{2})r)^{n+1}+((r^{2}+s^{2})r)^{n}}{2}=0 \\ 
\\ 
\dfrac{r^{n+1}-r^{n}}{\Delta t}+\dfrac{s_{xx}^{n+1}+s_{xx}^{n}}{2}+q\dfrac{%
((r^{2}+s^{2})s)^{n+1}+((r^{2}+s^{2})s)^{n}}{2}=0%
\end{array}
\label{CN}
\end{equation}%
where $\Delta t$ is the time step and superscripts denote the time levels.
For the numerical purpose, the nonlinear terms in system (\ref{CN}) can be
linearized by the technique in \cite{Rubin}.%
\begin{equation*}
\begin{array}{l}
((r^{2}+s^{2})r)^{n+1}=(r^{3})^{n+1}+(s^{2}r)^{n+1}=3\left( r^{n}\right)
^{2}r^{n+1}-2\left( r^{n}\right) ^{3}+2r^{n}s^{n}s^{n+1}+\left( s^{n}\right)
^{2}r^{n+1}-2\left( s^{n}\right) ^{2}r^{n}, \\ 
\\ 
((r^{2}+s^{2})s)^{n+1}=(r^{2}s)^{n+1}+(s^{3})^{n+1}=2s^{n}r^{n}r^{n+1}+%
\left( r^{n}\right) ^{2}s^{n+1}-2\left( r^{n}\right) ^{2}s^{n}+3\left(
s^{n}\right) ^{2}s^{n+1}-2\left( s^{n}\right) ^{3}.%
\end{array}%
\end{equation*}%
Substitution of the approximations of $r$ and $s$ with their related
derivatives into system (\ref{CN}) leads to the following fully discretized
equations%
\begin{equation}
\begin{array}{ll}
& \nu _{m1}\delta _{m-1}^{n+1}+\nu _{m2}\phi _{m-1}^{n+1}+\nu _{m3}\delta
_{m}^{n+1}+\nu _{m4}\phi _{m}^{n+1}+\nu _{m5}\delta _{m+1}^{n+1}+\nu
_{m6}\phi _{m+1}^{n+1} \\ 
= &  \\ 
& \nu _{m7}\delta _{m-1}^{n}+\nu _{m8}\phi _{m-1}^{n}+\nu _{m9}\delta
_{m}^{n}+\nu _{m10}\phi _{m}^{n}+\nu _{m11}\delta _{m+1}^{n}+\nu _{m12}\phi
_{m+1}^{n}%
\end{array}
\label{Sys_1}
\end{equation}%
and%
\begin{equation}
\begin{array}{ll}
& \nu _{m13}\delta _{m-1}^{n+1}+\nu _{m14}\phi _{m-1}^{n+1}+\nu _{m15}\delta
_{m}^{n+1}+\nu _{m16}\phi _{m}^{n+1}+\nu _{m17}\delta _{m+1}^{n+1}+\nu
_{m18}\phi _{m+1}^{n+1} \\ 
= &  \\ 
& \nu _{m19}\delta _{m-1}^{n}+\nu _{m20}\phi _{m-1}^{n}+\nu _{m21}\delta
_{m}^{n}+\nu _{m22}\phi _{m}^{n}+\nu _{m23}\delta _{m+1}^{n}+\nu _{m24}\phi
_{m+1}^{n}%
\end{array}
\label{Sys_2}
\end{equation}%
where%
\begin{equation*}
\begin{array}{lll}
\nu _{m1}=-\Delta t(q(3r^{2}+s^{2})\alpha _{1}+\gamma _{1}), &  & \nu
_{m13}=(2+2\Delta tqrs)\alpha _{1}, \\ 
\nu _{m2}=(2-2\Delta tqrs)\alpha _{1}, &  & \nu _{m14}=\Delta
t(q(r^{2}+s^{2})\alpha _{1}+\gamma _{1}), \\ 
\nu _{m3}=-\Delta t(q(3r^{2}+s^{2})\alpha _{2}+\gamma _{2}), &  & \nu
_{m15}=(2+2\Delta tqrs)\alpha _{2}, \\ 
\nu _{m4}=(2-2\Delta tqrs)\alpha _{2}, &  & \nu _{m16}=\Delta
t(q(r^{2}+s^{2})\alpha _{2}+\gamma _{2}), \\ 
\nu _{m5}=-\Delta t(q(3r^{2}+s^{2})\alpha _{1}+\gamma _{1}), &  & \nu
_{m17}=(2+2\Delta tqrs)\alpha _{1}, \\ 
\nu _{m6}=(2-2\Delta tqrs)\alpha _{1}, &  & \nu _{m18}=\Delta
t(q(r^{2}+s^{2})\alpha _{1}+\gamma _{1}), \\ 
\nu _{m7}=-\Delta tqr^{2}\alpha _{1}+\Delta t\gamma _{1}, &  & \nu
_{m19}=(2+\Delta tqrs)\alpha _{1}, \\ 
\nu _{m8}=(2-\Delta tqrs)\alpha _{1}, &  & \nu _{m20}=-\Delta tqs^{2}\alpha
_{1}-\Delta t\gamma _{1}, \\ 
\nu _{m9}=-\Delta tqr^{2}\alpha _{2}+\Delta t\gamma _{2}, &  & \nu
_{m21}=(2+\Delta tqrs)\alpha _{2}, \\ 
\nu _{m10}=(2-\Delta tqrs)\alpha _{2}, &  & \nu _{m22}=-\Delta tqs^{2}\alpha
_{2}-\Delta t\gamma _{2}, \\ 
\nu _{m11}=-\Delta tqr^{2}\alpha _{1}+\Delta t\gamma _{1}, &  & \nu
_{m23}=(2+\Delta tqrs)\alpha _{1}, \\ 
\nu _{m12}=(2-\Delta tqrs)\alpha _{1}, &  & \nu _{m24}=-\Delta tqs^{2}\alpha
_{1}-\Delta t\gamma _{1},%
\end{array}%
\end{equation*}%
\begin{equation*}
\begin{array}{cc}
r=\alpha _{1}\delta _{m-1}^{n}+\alpha _{2}\delta _{m}^{n}+\alpha _{3}\delta
_{m+1}^{n}, & s=\alpha _{1}\phi _{m-1}^{n}+\alpha _{2}\phi _{m}^{n}+\alpha
_{3}\phi _{m+1}^{n}%
\end{array}%
\end{equation*}%
and%
\begin{eqnarray*}
\alpha _{1} &=&\dfrac{s-ph}{2(phc-s)},\text{ \ \ \ \ \ \ }\alpha _{2}=1, \\
\gamma _{1} &=&\dfrac{p^{2}s}{2(phc-s)},\text{ \ \ \ \ \ \ }\gamma _{2}=-%
\dfrac{p^{2}s}{phc-s}.
\end{eqnarray*}%
There are $2N+2$ equations and $2N+6$ unknowns in systems (\ref{Sys_1}) and (%
\ref{Sys_2}). For the solvability of this system, the number of equations
and the number of unknown parameters should be equalized. The boundary
conditions enable us to eliminate the boundary parameters $\delta
_{-1}^{n+1},\delta _{N+1}^{n+1}$ and $\phi _{-1}^{n+1},$ $\phi _{N+1}^{n+1}$
from the system (\ref{Sys_1}) and (\ref{Sys_2}) such that we obtain a
solvable matrix system. To start the iteration, determination of the initial
parameters $\delta _{m}^{0}$ and $\phi _{m}^{0}$ are necessary. Once the
initial parameters are calculated then the time evolutions of the unknowns
are found from the recurrence relationship (\ref{Sys_1}) and (\ref{Sys_2}).

\section{Test problems}

This section is devoted for the observation of the efficiency of the method
so that several test problems are considered in order to illustrate the
accuracy. For this purpose we first calculate the possible $L_{\infty }$
error norm which is defined by%
\begin{equation*}
L_{\infty }=\underset{i}{\max }\left\vert
U_{i}^{exact}-U_{i}^{numerical}\right\vert ,
\end{equation*}

Additionally, the invariants of Eq.(\ref{NLS}) also give an idea about the
accuracy of the method especially in cases that the equation does not have
an analytical solution. Although there are infinitely many conservation laws
for Eq.(\ref{NLS}), here we investigate only the followings: 
\begin{equation*}
\begin{tabular}{ll}
C$_{1}=\int\limits_{a}^{b}\left\vert U\right\vert ^{2}dx,$ & C$%
_{2}=\int\limits_{a}^{b}\left( \left\vert U_{x}\right\vert ^{2}-\dfrac{1}{2}%
q\left\vert U\right\vert ^{4}\right) dx,$%
\end{tabular}%
\end{equation*}

\subsection{Single soliton}

The function%
\begin{equation}
U(x,t)=\alpha \sqrt{2/q}e^{i\left( \frac{S}{2}x-\frac{1}{4}(S^{2}-\alpha
^{2})t\right) }\text{sech}\left( \alpha (x-St)\right)   \label{Analytic}
\end{equation}%
represents the single soliton solution of Eq.(\ref{NLS}). When $t$ is fixed
then the solution (\ref{Analytic}) decays exponentially as $\left\vert
x\right\vert \rightarrow \infty $. The initial and the boundary conditions
are inferenced from the above solution. 
\begin{equation*}
\begin{tabular}{cc}
\FRAME{itbpF}{3.0441in}{2.5071in}{0in}{}{}{fig2.eps}{\special{language
"Scientific Word";type "GRAPHIC";maintain-aspect-ratio TRUE;display
"USEDEF";valid_file "F";width 3.0441in;height 2.5071in;depth
0in;original-width 3in;original-height 2.4664in;cropleft "0";croptop
"1";cropright "1";cropbottom "0";filename 'Fig2.eps';file-properties
"XNPEU";}} & \FRAME{itbpF}{3.0441in}{2.5227in}{0in}{}{}{fig3.eps}{\special%
{language "Scientific Word";type "GRAPHIC";maintain-aspect-ratio
TRUE;display "USEDEF";valid_file "F";width 3.0441in;height 2.5227in;depth
0in;original-width 3in;original-height 2.4811in;cropleft "0";croptop
"1";cropright "1";cropbottom "0";filename 'Fig3.eps';file-properties
"XNPEU";}} \\ 
\textbf{Fig.2: }Single soliton profiles & \textbf{Fig.3: }Error distribution
at for $p=0.0000182$%
\end{tabular}%
\end{equation*}

Eq.(\ref{Analytic}) gives a wave that moves with speed $S$ and its magnitude
is governed by the real parameter $\alpha $. Since it is a useful tool for
comparison, this problem is a well known example in the literature. Table 2
presents a detailed comparison on this example for different parameter
choices. The solution profiles and the absolute error distribution are
illustrated in Figs.2-3.%
\begin{equation*}
\begin{tabular}{|l|l|l|l|l|}
\hline
\multicolumn{5}{|l|}{\textbf{Table 2}} \\ \hline
\multicolumn{5}{|l|}{Errors at $t=1$ for $q=2,$ $S=4,$ $\alpha =1$} \\ \hline
Method & $h$ & $\Delta t$ & $L_{\infty }$ ($p=1$) & $L_{\infty }$ (various $p
$) \\ \hline
Present & $0,05$ & $0,005$ & $0,0057$ & $0,0015104$ \\ \hline
&  &  &  & ($p=0,0000078348$) \\ \hline
& $0,3125$ & $0,02$ & $0,1872$ & $0,0064891$ \\ \hline
&  &  &  & ($p=0,0000001289$) \\ \hline
& $0,3125$ & $0,0026$ & $0,1913$ & $0,0053$ \\ \hline
&  &  &  & ($p=0,0000001289$) \\ \hline
& $0,06$ & $0,0165$ & $0,0048$ & $0.0015$ \\ \hline
&  &  &  & ($p=0,0000031976$) \\ \hline
& $0,05$ & $0,04$ & $0,0168$ & $0.0039$ \\ \hline
&  &  &  & ($p=0.0000020600$) \\ \hline
Kuintik B-Spline (Saka, 2012) & $0,05$ & $0,005$ &  & $0,0003$ \\ \hline
& $0,3125$ & $0,02$ &  & $0,002$ \\ \hline
& $0,3125$ & $0,0026$ &  & $0,006$ \\ \hline
B-spline Galerkin (Da\u{g}, 1999) & $0,05$ & $0,005$ &  & $0,0003$ \\ \hline
& $0,3125$ & $0,02$ &  & $0,002$ \\ \hline
B-spline Col. (Gardner vd., 1993) & $0,05$ & $0,005$ &  & $0,008$ \\ \hline
& $0,03$ & $0,005$ &  & $0,002$ \\ \hline
Kapal\i\ (C-N) (Taha ve Ablowitz, 1984) & $0,05$ & $0,005$ &  & $0,00585$ \\ 
\hline
Split step Fourier (Taha ve Ablowitz, 1984) & $0,3125$ & $0,02$ &  & $0,00466
$ \\ \hline
A-L local (Taha ve Ablowitz, 1984) & $0,06$ & $0,0165$ &  & $0,00580$ \\ 
\hline
A-L global (Taha ve Ablowitz, 1984) & $0,05$ & $0,04$ &  & $0,00561$ \\ 
\hline
Pseudospectral (Taha ve Ablowitz, 1984) & $0,3125$ & $0,0026$ &  & $0,00513$
\\ \hline
\end{tabular}%
\end{equation*}

\subsection{Collision of two solitons that move in opposite directions}

A collision of two solitons that travels in opposite directions can be
observed with the initial conditions%
\begin{equation}
U(x,0)=U_{1}(x,0)+U_{2}(x,0),  \label{Collision}
\end{equation}%
where%
\begin{equation}
U_{j}(x,0)=\alpha _{j}\sqrt{2/q}e^{i\left( \frac{S}{2}(x-x_{j})\right) }%
\text{sech}\left( \alpha _{j}(x-x_{j})\right) ,\text{ }j=1,2.
\label{Collision_Initials}
\end{equation}%
The parameters in this test problem are considered as%
\begin{equation*}
q=2,h=0.1,\Delta t=0.005,\alpha _{1}=1.0,S_{1}=-4.0,x_{1}=10,\alpha
_{2}=1.0,S_{2}=4.0,x_{2}=-10
\end{equation*}%
to coincide with those of some earlier studies. Since there is no available
analytical soliton satisfying the given initial condition, only the
invariants are considered for the observation of the accuracy. According to
the results shown in Table 3, it can be concluded that the method produces
acceptable results.

\bigskip 
\begin{equation*}
\begin{tabular}{lcllll}
\multicolumn{6}{l}{\textbf{Table 3}} \\ 
\multicolumn{6}{l}{Invariants for the collision of two solitons} \\ \hline
& Time &  & C$_{1}$ &  & C$_{2}$ \\ \hline
$p=1$ & $0.0$ &  & \multicolumn{1}{c}{$3.99999$} &  & \multicolumn{1}{c}{$%
14.66577$} \\ 
& $0.5$ &  & \multicolumn{1}{c}{$4.00000$} &  & \multicolumn{1}{c}{$14.66634$%
} \\ 
& $1.0$ &  & \multicolumn{1}{c}{$4.00000$} &  & \multicolumn{1}{c}{$14.66706$%
} \\ 
& $1.5$ &  & \multicolumn{1}{c}{$4.00001$} &  & \multicolumn{1}{c}{$14.66761$%
} \\ 
& $2.0$ &  & \multicolumn{1}{c}{$4.00001$} &  & \multicolumn{1}{c}{$14.66694$%
} \\ 
& $2.5$ &  & \multicolumn{1}{c}{$3.99992$} &  & \multicolumn{1}{c}{$14.61149$%
} \\ 
& $3.0$ &  & \multicolumn{1}{c}{$4.00003$} &  & \multicolumn{1}{c}{$14.66809$%
} \\ 
& $3.5$ &  & \multicolumn{1}{c}{$4.00003$} &  & \multicolumn{1}{c}{$14.66803$%
} \\ 
& $4.0$ &  & \multicolumn{1}{c}{$4.00004$} &  & \multicolumn{1}{c}{$14.66766$%
} \\ 
& $4.5$ &  & \multicolumn{1}{c}{$4.00005$} &  & \multicolumn{1}{c}{$14.66734$%
} \\ 
& $5.0$ &  & \multicolumn{1}{c}{$4.00005$} &  & \multicolumn{1}{c}{$14.66705$%
} \\ 
& $5.5$ &  & \multicolumn{1}{c}{$4.00006$} &  & \multicolumn{1}{c}{$14.66691$%
} \\ 
& $6.0$ &  & \multicolumn{1}{c}{$4.00007$} &  & \multicolumn{1}{c}{$14.66690$%
} \\ \hline
\cite{Aksoy} $\lambda =0$ & $6.0$ &  & $4.00000$ &  & $14.66669$ \\ \hline
\end{tabular}%
\end{equation*}

\bigskip

Eq.(\ref{Collision_Initials}) represents two solitons having equal
magnitudes and velocities that are $1$ and $4$ respectively. One of these
waves is located at $x=-10$ whereas the other is at $x=10$. The waves move
in opposite directions starting from their mentioned initial locations.
Fig.4 illustrates the before and the later positions of interacting waves.
On the other hand the collision is monitored in Fig.5. As an expected
situation, these figures show that the solitons keep their initial profile
and properites after the collision.

\bigskip 
\begin{equation*}
\begin{tabular}{ll}
\FRAME{itbpF}{3.0441in}{2.578in}{0in}{}{}{fig4.eps}{\special{language
"Scientific Word";type "GRAPHIC";maintain-aspect-ratio TRUE;display
"USEDEF";valid_file "F";width 3.0441in;height 2.578in;depth
0in;original-width 3in;original-height 2.5365in;cropleft "0";croptop
"1";cropright "1";cropbottom "0";filename 'Fig4.eps';file-properties
"XNPEU";}} & \FRAME{itbpF}{3.0441in}{2.6057in}{0in}{}{}{fig5.eps}{\special%
{language "Scientific Word";type "GRAPHIC";maintain-aspect-ratio
TRUE;display "USEDEF";valid_file "F";width 3.0441in;height 2.6057in;depth
0in;original-width 3in;original-height 2.5642in;cropleft "0";croptop
"1";cropright "1";cropbottom "0";filename 'Fig5.eps';file-properties
"XNPEU";}} \\ 
\multicolumn{1}{c}{\textbf{Fig.4:} Waves at before and later the collision}
& \multicolumn{1}{c}{\textbf{Fig.5:} Wave profiles at just before the
collision}%
\end{tabular}%
\end{equation*}

\bigskip

\begin{equation*}
\begin{tabular}{l}
\FRAME{itbpF}{3.0441in}{2.5374in}{0in}{}{}{fig6.eps}{\special{language
"Scientific Word";type "GRAPHIC";maintain-aspect-ratio TRUE;display
"USEDEF";valid_file "F";width 3.0441in;height 2.5374in;depth
0in;original-width 3in;original-height 2.4967in;cropleft "0";croptop
"1";cropright "1";cropbottom "0";filename 'Fig6.eps';file-properties
"XNPEU";}} \\ 
\multicolumn{1}{c}{\textbf{Fig.6:} Wave profiles at interaction}%
\end{tabular}%
\end{equation*}

\subsection{Birth of soliton}

Here we focus on two different examples which are the birth of standing
soliton and the birth of mobile soliton.

\subsubsection{Birth of standing soliton}

According to the theory, an initial condition satisfying%
\begin{equation*}
C=\dint\limits_{-\infty }^{\infty }U(x,0)dx\geq \pi ,
\end{equation*}%
results a soliton in time process. Otherwise the soliton decays away. The
verification can be seen for some different numerical methods in \cite%
{gr,idag,Saka}. This observation can be simulated with the Maxwellian
initial condition%
\begin{equation}
U(x,0)=Ae^{-x^{2}}  \label{Maxwellian1}
\end{equation}%
or $C=\sqrt{\pi }A$, so the usege of (\ref{Maxwellian1}) with $A\geq \sqrt{%
\pi }\approx 1.7725$ produces a soliton whereas the choice of $A<\sqrt{\pi }$
yields a fading out initial condition. Figs.7 and 8 show this procedure. The
invariants are listed in Table 4 for comparison. The analytical invariants
are 
\begin{eqnarray*}
C_{1} &=&A^{2}\sqrt{\frac{\pi }{2}}=3.9710, \\
C_{2} &=&\frac{1}{4}A^{2}(2\sqrt{2}-qA^{2})\sqrt{\pi }=-4.9256,
\end{eqnarray*}%
so it can be concluded that the presented results are in good agreement with
the exact invariants.

\bigskip 
\begin{equation*}
\begin{tabular}{cc}
\FRAME{itbpF}{3.3883in}{3.0606in}{0in}{}{}{fig7.bmp}{\special{language
"Scientific Word";type "GRAPHIC";maintain-aspect-ratio TRUE;display
"USEDEF";valid_file "F";width 3.3883in;height 3.0606in;depth
0in;original-width 3.3607in;original-height 3.0329in;cropleft "0";croptop
"1";cropright "1";cropbottom "0";filename 'Fig7.bmp';file-properties
"XNPEU";}} & \FRAME{itbpF}{3.3814in}{3.1073in}{0in}{}{}{fig8.bmp}{\special%
{language "Scientific Word";type "GRAPHIC";maintain-aspect-ratio
TRUE;display "USEDEF";valid_file "F";width 3.3814in;height 3.1073in;depth
0in;original-width 3.3529in;original-height 3.0796in;cropleft "0";croptop
"1";cropright "1";cropbottom "0";filename 'Fig8.bmp';file-properties
"XNPEU";}} \\ 
\textbf{Fig7: }Formation of standing solution for $A=1$ & \textbf{Fig8:}%
Formation of standing solution for $A=1.78$%
\end{tabular}%
\end{equation*}

\bigskip

\begin{equation*}
\begin{tabular}{llllcllllllllll}
\multicolumn{15}{l}{\textbf{Table 4}} \\ 
\multicolumn{15}{l}{Invariants for $A=1.78$ and $-45\leq x\leq 45,$ $N=1334,$
$\Delta t=0.005,q=2.$} \\ \hline
&  &  &  & Time &  &  &  &  & C$_{1}$ &  &  &  &  & C$_{2}$ \\ \hline
$p=1$ &  &  &  & $0.0$ &  &  &  &  & $3.97100$ &  &  &  &  & $-4.92563$ \\ 
&  &  &  & $2.0$ &  &  &  &  & $3.97091$ &  &  &  &  & $-4.93275$ \\ 
&  &  &  & $4.0$ &  &  &  &  & $3.97076$ &  &  &  &  & $-4.93191$ \\ 
&  &  &  & $6.0$ &  &  &  &  & $3.97062$ &  &  &  &  & $-4.93128$ \\ \hline
\cite{Aksoy} $\lambda =0$ &  &  &  & $6.0$ &  &  &  &  & $3.97093$ &  &  & 
&  & $-4.92672$ \\ \hline
\end{tabular}%
\end{equation*}

\subsubsection{Birth of mobile soliton}

A mobile soliton can be studied with the Maxwellian initial condition%
\begin{equation}
U(x,0)=Ae^{-x^{2}+2ix}  \label{Maxwellian2}
\end{equation}%
which produces a mobile soliton having the velocity $4$ and height $2$ with
parameter choice $A=1.78.$ The same parameter choices with the previous
problem form a soliton that has the peak position at $x=24\ $when $t=6.$ The
travelling soliton wave is graphed in Fig.9 for $A=1.78.$ As seen from
Fig.10 that the case $A=1$ does not produce any soliton. 
\begin{equation*}
\begin{tabular}{ll}
\FRAME{itbpF}{3.3883in}{3.0139in}{0in}{}{}{fig9.bmp}{\special{language
"Scientific Word";type "GRAPHIC";maintain-aspect-ratio TRUE;display
"USEDEF";valid_file "F";width 3.3883in;height 3.0139in;depth
0in;original-width 3.3607in;original-height 2.9871in;cropleft "0";croptop
"1";cropright "1";cropbottom "0";filename 'Fig9.bmp';file-properties
"XNPEU";}} & \FRAME{itbpF}{3.3079in}{3.0874in}{0in}{}{}{fig10.bmp}{\special%
{language "Scientific Word";type "GRAPHIC";maintain-aspect-ratio
TRUE;display "USEDEF";valid_file "F";width 3.3079in;height 3.0874in;depth
0in;original-width 3.2802in;original-height 3.0597in;cropleft "0";croptop
"1";cropright "1";cropbottom "0";filename 'Fig10.bmp';file-properties
"XNPEU";}} \\ 
\textbf{Fig9: }Formation of travelling solution for$A=1$ & \textbf{Fig10: }%
Formation of travelling solution for $A=1.78$%
\end{tabular}%
\end{equation*}

\bigskip

Analytical invariants for this case are calculated as%
\begin{eqnarray*}
C_{1} &=&\sqrt{\dfrac{\pi }{2}}A^{2}=3.97100, \\
C_{2} &=&5\sqrt{\dfrac{\pi }{2}}A^{2}-\dfrac{\sqrt{\pi }}{4}qA^{4}=10.95838.
\end{eqnarray*}

Table 5 indicates that $C_{1}$ has been found to be constant and there is
variation in $C_{2}$ less than $1.95\%.$

\bigskip

\begin{equation*}
\begin{tabular}{llllcllllllllll}
\multicolumn{15}{l}{\textbf{Table 5}} \\ 
\multicolumn{15}{l}{Invariants for $A=1.78$ and $-45\leq x\leq 45,$ $N=1334,$
$\Delta t=0.005,$ $q=2.$} \\ \hline
&  &  &  & Time &  &  &  &  & C$_{1}$ &  &  &  &  & C$_{2}$ \\ \hline
$p=1$ &  &  &  & $0.0$ &  &  &  &  & $3.97100$ &  &  &  &  & $10.95788$ \\ 
&  &  &  & $2.0$ &  &  &  &  & $3.97113$ &  &  &  &  & $10.94128$ \\ 
&  &  &  & $4.0$ &  &  &  &  & $3.97112$ &  &  &  &  & $10.94298$ \\ 
&  &  &  & $6.0$ &  &  &  &  & $3.97113$ &  &  &  &  & $10.94343$ \\ \hline
\cite{Aksoy} $\lambda =0$ &  &  &  & $4.0$ &  &  &  &  & $3.97087$ &  &  & 
&  & $10.95432$ \\ \hline
\end{tabular}%
\end{equation*}

\bigskip

\subsection{Bound state solution}

It is stated in \cite{gar} that if the condition%
\begin{equation*}
q=2M^{2},\text{ }M=1,2,...
\end{equation*}%
holds then the initial condition%
\begin{equation*}
U\left( x,0\right) =\sec \text{h}\left( x\right)
\end{equation*}%
evolves $M$ soliton waves. Illustrations of these solitons are depicted in
Figs.11-18 where the each wave profile and their trajectories, to see the
density variations, are given alongside. The numerical computations have
been carried out for%
\begin{equation*}
h=0.03,\Delta t=0.005,N=1334,q=32,50,-20\leq x\leq 20.
\end{equation*}

\bigskip 
\begin{equation*}
\begin{tabular}{cc}
\FRAME{itbpF}{3.7075in}{3.7075in}{0in}{}{}{fig11.bmp}{\special{language
"Scientific Word";type "GRAPHIC";maintain-aspect-ratio TRUE;display
"USEDEF";valid_file "F";width 3.7075in;height 3.7075in;depth
0in;original-width 3.6599in;original-height 3.6599in;cropleft "0";croptop
"1";cropright "1";cropbottom "0";filename 'Fig11.bmp';file-properties
"XNPEU";}} & \FRAME{itbpF}{3.7075in}{3.7075in}{0in}{}{}{fig12.bmp}{\special%
{language "Scientific Word";type "GRAPHIC";maintain-aspect-ratio
TRUE;display "USEDEF";valid_file "F";width 3.7075in;height 3.7075in;depth
0in;original-width 3.6599in;original-height 3.6599in;cropleft "0";croptop
"1";cropright "1";cropbottom "0";filename 'Fig12.bmp';file-properties
"XNPEU";}} \\ 
\textbf{Fig.11:} Wave profiles for $M=4$ & \textbf{Fig.12:} Trajectories for 
$M=4$%
\end{tabular}%
\end{equation*}

\bigskip 
\begin{equation*}
\begin{tabular}{cc}
\FRAME{itbpF}{3.7075in}{3.7075in}{0in}{}{}{fig13.bmp}{\special{language
"Scientific Word";type "GRAPHIC";maintain-aspect-ratio TRUE;display
"USEDEF";valid_file "F";width 3.7075in;height 3.7075in;depth
0in;original-width 3.6599in;original-height 3.6599in;cropleft "0";croptop
"1";cropright "1";cropbottom "0";filename 'Fig13.bmp';file-properties
"XNPEU";}} & \FRAME{itbpF}{3.7075in}{3.7075in}{0in}{}{}{fig14.bmp}{\special%
{language "Scientific Word";type "GRAPHIC";maintain-aspect-ratio
TRUE;display "USEDEF";valid_file "F";width 3.7075in;height 3.7075in;depth
0in;original-width 3.6599in;original-height 3.6599in;cropleft "0";croptop
"1";cropright "1";cropbottom "0";filename 'Fig14.bmp';file-properties
"XNPEU";}} \\ 
\textbf{Fig.13:} Wave profiles for $M=5$ & \textbf{Fig.14:} Trajectories for 
$M=5$%
\end{tabular}%
\end{equation*}%
\begin{equation*}
\begin{tabular}{cc}
\FRAME{itbpF}{3.7075in}{3.7075in}{0in}{}{}{fig15.bmp}{\special{language
"Scientific Word";type "GRAPHIC";maintain-aspect-ratio TRUE;display
"USEDEF";valid_file "F";width 3.7075in;height 3.7075in;depth
0in;original-width 3.6599in;original-height 3.6599in;cropleft "0";croptop
"1";cropright "1";cropbottom "0";filename 'Fig15.bmp';file-properties
"XNPEU";}} & \FRAME{itbpF}{3.7075in}{3.7075in}{0in}{}{}{fig16.bmp}{\special%
{language "Scientific Word";type "GRAPHIC";maintain-aspect-ratio
TRUE;display "USEDEF";valid_file "F";width 3.7075in;height 3.7075in;depth
0in;original-width 3.6599in;original-height 3.6599in;cropleft "0";croptop
"1";cropright "1";cropbottom "0";filename 'Fig16.bmp';file-properties
"XNPEU";}} \\ 
\textbf{Fig.15:} Wave profiles for $M=6$ & \textbf{Fig.16:} Trajectories for 
$M=6$%
\end{tabular}%
\end{equation*}

\bigskip 
\begin{equation*}
\begin{tabular}{cc}
\FRAME{itbpF}{3.7075in}{3.7075in}{0in}{}{}{fig17.bmp}{\special{language
"Scientific Word";type "GRAPHIC";maintain-aspect-ratio TRUE;display
"USEDEF";valid_file "F";width 3.7075in;height 3.7075in;depth
0in;original-width 3.6599in;original-height 3.6599in;cropleft "0";croptop
"1";cropright "1";cropbottom "0";filename 'Fig17.bmp';file-properties
"XNPEU";}} & \FRAME{itbpF}{3.7075in}{3.7075in}{0in}{}{}{fig18.bmp}{\special%
{language "Scientific Word";type "GRAPHIC";maintain-aspect-ratio
TRUE;display "USEDEF";valid_file "F";width 3.7075in;height 3.7075in;depth
0in;original-width 3.6599in;original-height 3.6599in;cropleft "0";croptop
"1";cropright "1";cropbottom "0";filename 'Fig18.bmp';file-properties
"XNPEU";}} \\ 
\textbf{Fig.17:} Wave profiles for $M=7$ & \textbf{Fig.18:} Trajectories for 
$M=7$%
\end{tabular}%
\end{equation*}

\bigskip

\section{Conclusion}

In this paper, exponential cubic B-spline collocation method is implemented
in order to get the solution of the cubic nonlinear Schr\"{o}dinger
equation. Over the uniform mesh, Crank-Nicolson formulas are employed for
time discretization whereas Rubin and Graves\cite{Rubin} technique is used
for the linearization. Four test problems that related to single soliton
wave, interaction of two opposite solitons, birth of soliton and the bound
state solution are examined for testing the numerical scheme. Comparisons
between the obtained results and some earlier papers show that the present
results are all acceptable and in agreement with those in the literature.
Simple adaptation and yielding band matrices can be stated as the advantages
of the method. On the other hand, requiring the determination of the free
parameter $p$ is an undesirable situation. In conclusion, exponential cubic
B-spline collocation method can be considered as a conservative numerical
method that leads to reasonable results

\end{document}